\documentclass[a4paper,12pt,reqno]{amsart}
\usepackage[headings]{fullpage}

\usepackage{hyperref}

\usepackage{amsthm}
\usepackage{amsmath}
\usepackage{amssymb}

\usepackage[all]{xy}   

\theoremstyle{plain}
\newtheorem{thm}{Theorem}[section]

\newtheorem{cor}[thm]{Corollary}

\theoremstyle{definition}
\newtheorem{defn}[thm]{Definition}

\theoremstyle{remark}
\newtheorem{rem}[thm]{Remark}

\newtheorem*{ackn*}{Acknowledgements}

\newcommand{\Hom}{\operatorname{Hom}}
\newcommand{\Vol}{\operatorname{Vol}}

\begin{document}

\title{On hyperbolic cohomology classes}

\author{M.~Brunnbauer}
\address{Mathematisches Institut, {\smaller LMU} M\"unchen, Theresienstr.~39, 80333~M\"unchen, Germany}
\email{michael.brunnbauer@mathematik.uni-muenchen.de}

\author{D.~Kotschick}
\address{Mathematisches Institut, {\smaller LMU} M\"unchen, Theresienstr.~39, 80333~M\"unchen, Germany}
\email{dieter@member.ams.org}

\date{\today; \copyright{\ M.~Brunnbauer and D.~Kotschick 2008}}
\thanks{The authors are members of the DFG priority program in Global Differential Geometry.}

%\keywords{}
%Primary; Secondary}

\begin{abstract}
We study hyperbolic cohomology classes in the general context of simplicial complexes and prove homological invariance statements for them. 
We relate the existence of hyperbolic cohomology classes to the non-amenability of the fundamental group. In degree two we clarify the relation 
between hyperbolic and atoroidal classes, leading to an application to symplectically atoroidal manifolds.
\end{abstract}

\maketitle

%%%%%%%%%%%%%%%%%%%%%%%%%%%%%%%%%%%%%%%%%%%%%%%%%%%%%%%%%%%%%%%%%%%%%%%%%%%%%%%%%%%%%
%%%%%%%%%%%%%%%%%%%%%%%%%%%%%%%%%%%%%%%%%%%%%%%%%%%%%%%%%%%%%%%%%%%%%%%%%%%%%%%%%%%%%

\section{Introduction}\label{s:intro}

The notion of hyperbolic cohomology classes was introduced by Gromov~\cite{Gromov(1991)} under the name of $\tilde d(\textrm{bounded})$ classes. The paper~\cite{Gromov(1991)} was concerned with K\"ahler manifolds whose K\"ahler classes are hyperbolic, and the more general case of hyperbolic classes represented by symplectic forms which need not be K\"ahler was considered by Polterovich~\cite{Pol02} and, more recently, by K{\c{e}}dra~\cite{Kedra(2007)}. The present note was motivated by these papers, and it aims to answer some of the questions raised, explicitly or implicitly, in~\cite{Kedra(2007)}.

In Section~\ref{s:classes} we investigate hyperbolic classes in the general context of simplicial complexes, without restricting to smooth manifolds. Our discussion proceeds along the lines of Gromov's ideas as presented in~\cite{Gromov(1991),Gromov(1993)}. In Theorems~\ref{t:hominv} (for degree two) and \ref{t:higherhominv} (for higher degrees) we prove a kind of homological invariance with respect to classifying maps for the notion of hyperbolicity of a cohomology class. These results, which follow a pattern of results on homological invariance of other largeness properties established in the work of the  first author~\cite{Brunnbauer(2007),BrunnbauerThesis}, clarify the discussions in~\cite[Section 0.2]{Gromov(1991)} and in~\cite[Section 5]{Kedra(2007)}.

In Section~\ref{amenable} we discuss the relation between amenability of the fundamental group and the absence of hyperbolic cohomology classes using ideas of Brooks~\cite{Brooks(1983),Brooks(1985)}. 
%and Gromov~\cite{Gromov(1991),Gromov(1993)}. 
Finally in Section~\ref{s:examples} we give concrete answers to some questions formulated in~\cite{Kedra(2007)}.

\section{Aspherical, atoroidal, hyperbolic, and bounded classes}\label{s:classes}

Let $X$ be a topological space. A cohomology class $w\in H^k(X;\mathbb{R})$ is called \emph{aspherical} if its pullback to any sphere is zero. Using the universal coefficient theorem $H^k(X;\mathbb{R})=\Hom(H_k(X;\mathbb{R}),\mathbb{R})$, one sees that a cohomology class $w$ is aspherical if and only if it maps every homology class represented by a sphere to zero, that is, if the image of the Hurewicz homorphism $\pi_k(X)\to H_k(X;\mathbb{R})$ lies in kernel of $w$. The subspace of all aspherical cohomology classes in $H^k(X;\mathbb{R})$ will be denoted by $V^k_{asph}(X)$.

In degree two, we will also consider the subspace $V^2_{ator}(X)\subset H^2(X;\mathbb{R})$ of all \emph{atoroidal} cohomology classes. A class is called atoroidal if it evaluates to zero on every homology class represented by a $2$-torus.

Since there is a degree one map $T^2 \to S^2$ it follows immediately that
\[ V^2_{ator}(X) \subset V^2_{asph}(X). \]
Moreover, these subspaces are natural in the following sense: if $f:X\to Y$ is a continuous map, then $f^*V^k_{asph}(Y) \subset V^k_{asph}(X)$ and $f^*V^2_{ator}(Y) \subset V^2_{ator}(X)$.

Note that $V^1_{asph}(X)=0$ for every space $X$ since every integral $1$-cycle is represented by a loop, and that $V^k_{asph}(X)=H^k(X;\mathbb{R})$ for $k\geq 2$ if $X$ is aspherical (i.\,e.\ $\pi_k(X)=0$ for $k\geq 2$). Furthermore, tori provide examples of aspherical spaces for which the atoroidal subspace is trivial.

Next, we want to define the notion of hyperbolic cohomology class. This is only possible on spaces for which cohomology classes can be represented by differential forms. Beyond smooth manifolds, this works for simplicial complexes, see for example \cite{Swan(1975)}. 

Let $X$ be a simplicial complex. A $k$-form $\omega$ on $X$ consists of a smooth $k$-form $\omega_\sigma$ for every simplex $\sigma\subset X$ such that $\omega_\sigma |_\tau \equiv \omega_\tau$ whenever $\tau\subset\sigma$ is a subsimplex. The space of all $k$-forms is denoted by $\Omega^k(X)$. Together with the exterior derivative this defines the de\,Rham complex of $X$. The following is the de\,Rham theorem for this general context:
\begin{thm}[de\,Rham, Thom]
There is a natural isomorphism
\[ H^k_{dR} (X) \xrightarrow{\cong} H^k(X;\mathbb{R}) \]
between the cohomology of the de\,Rham complex and real singular cohomology.
\end{thm}

Analogously to the definition of differential forms, one defines Riemannian metrics for simplicial complexes. A Riemannian metric $g$ consists of a Riemannian metric $g_\sigma$ on every simplex $\sigma$ of $X$ such that $g_\sigma|_\tau \equiv g_\tau$ for $\tau\subset\sigma$. Using this, we can make the following definition following Gromov~\cite{Gromov(1991),Gromov(1993)}:
\begin{defn}\label{d:hyp_subsp}
Let $X$ be a finite simplicial complex. Denote the universal covering by $p\colon\tilde X\to X$. A cohomology class $w\in H^k(X;\mathbb{R})$ is called \emph{hyperbolic} if the pullback $p^*\omega$ of a representing $k$-form $\omega$ has a primitive $\alpha\in\Omega^{k-1}(\tilde X)$ which is bounded with respect to some lifted Riemannian metric.
\end{defn}
Note that this definition does not depend on the choice of the Riemannian metric since $X$ is compact. If the pullback of one representative has a bounded primitive, then this holds for every representative. (The difference of two representing forms is the exterior derivative of a form on $X$, whose lift to $\tilde X$ is obviously bounded.)

The subset $V^k_{hyp}(X)\subset H^k(X;\mathbb{R})$ of all hyperbolic classes is a subspace. If $f:X\to Y$ is a simplicial map, then obviously $f^*V^k_{hyp}(Y) \subset V^k_{hyp}(X)$. Since every continuous map may be approximated by a simplicial one, it follows that the hyperbolic subspace is natural with respect to continuous maps.

Note that $V^k_{hyp}(X)\subset V^k_{asph}(X)$ for $k\geq 2$ since every map $S^k\to X$ factorizes through the universal covering $\tilde X\to X$ and hyperbolic classes are by definition cohomologous to zero on $\tilde X$. In \cite{Kedra(2007)}, Proposition 1.9, K\c{e}dra showed that every hyperbolic class of degree two is atoroidal (see also Section~\ref{amenable} below). Moreover, note that $V^1_{hyp}(X)=0$. (This is a direct consequence of Theorem~\ref{amenable->hyp=0}, but it is also rather obvious.)

We will also consider bounded cohomology classes. Let $X$ be a topological space. Denote by $S_k(X)$ the set of all singular $k$-simplices in $X$. The space of singular $k$-cochains is given by $C^k(X;\mathbb{R})=\{ c:S_k(X)\to\mathbb{R} \}$, the vector space of all real functions on $S_k(X)$. The bounded cochain group $C^k_b(X)$ consists of all such functions which are uniformly bounded on $S_k(X)$. The bounded cohomology $H^k_b(X)$ is the cohomology of this subcomplex of the singular cochain complex. (More details and deep results on bounded cohomology can be found in Gromov's paper~\cite{Gromov(1983)}.)

The image of the canonical homomorphism $H^k_b(X)\to H^k(X;\mathbb{R})$ is denoted by $V^k_b(X)$. The cohomology classes in this subspace are called \emph{bounded}. This subspace is natural with respect to continuous maps, and K\c{e}dra~\cite[Theorem 2.1]{Kedra(2007)} proved that it is contained in the hyperbolic subspace. 
In fact, K\c{e}dra stated this only for closed manifolds $M$ in place of finite simplicial complexes $X$, but the result is true in this generality with the same proof. Since the bounded cohomology of spheres and of tori is trivial, it is clear that bounded classes are always atoroidal and aspherical.

Let $X$ be a finite simplicial complex. We have defined four subspaces of $H^k(X;\mathbb{R})$ and seen that they fulfill the following relations:
\begin{align*}
&V^k_b(X) \subset V^k_{hyp}(X) \subset V^k_{asph}(X) \;\;\;\;\;\mathrm{respectively} \\
V^2_b(X) &\subset V^2_{hyp}(X) \subset V^2_{ator}(X) \subset V^2_{asph}(X).
\end{align*}

Denote by $\pi$ the fundamental group of $X$ and by $c\colon X\to B\pi$ the classifying map of the universal covering. This is a map that induces the identity on fundamental groups and is uniquely determined up to homotopy by this condition. Without loss of generality we may assume that $c\colon X\hookrightarrow B\pi$ is the inclusion of a subcomplex such that the $2$-skeleton of $B\pi$ is contained in $X$. By the long exact cohomology sequence and the Hurewicz theorem, it follows that the induced homomorphism
\[ c^* \colon H^2 (B\pi;\mathbb{R}) \hookrightarrow H^2(X;\mathbb{R}) \]
is injective and the image of $c^*$ is $V^2_{asph}(X)$.

Since $V^2_{asph}(B\pi)=H^2(B\pi;\mathbb{R})$, this may be rephrased by $c^*(V^2_{asph}(B\pi)) = V^2_{asph}(X)$. By cellular approximation every map $T^2\to B\pi$ may be homotoped to a map $T^2\to X\subset B\pi$. Thus, if $c^*w$ is atoroidal, then $w$ has to be atoroidal too, that is, $c^*(V^2_{ator}(B\pi)) = V^2_{ator}(X)$. Furthermore, $c$ induces an isomorphism on bounded cohomology (\cite{Gromov(1983)}, page 40). Hence, $c^*(V^2_b(B\pi)) = V^2_b(X)$ follows.

In general, $B\pi$ does not have the homotopy type of a finite complex (for example if $\pi$ contains non-trivial torsion elements). Therefore, the hyperbolic subspace is not defined for $B\pi$ by Definition \ref{d:hyp_subsp}. That definition does not apply, because when a simplicial complex is not finite, different metrics are not bi-Lipschitz equivalent, and so the notion of bounded primitive depends on the choice of metric. Nevertheless, there is the following definition due to Gromov~\cite[Subsection 0.2.C]{Gromov(1991)}:
\begin{defn}
A cohomology class $w\in H^k(B\pi;\mathbb{Q})$ is called \emph{hyperbolic} if its pullback to any finite simplicial complex is hyperbolic.
\end{defn}
The next theorem shows that in the situation above where $X\subset B\pi$ contains the $2$-skeleton the equality
\begin{equation}\label{eq:hom_invar_hyp}
c^* V^2_{hyp}(B\pi) = V^2_{hyp}(X)
\end{equation}
holds. Thus, the four subspaces in two-dimensional cohomology depend only on the fundamental group and on the classifying map of the universal covering.

\begin{thm}\label{t:hominv}
Let $X$ and $Y$ be two finite simplicial complexes, let $c:X\to B\pi$ be the classifying map of the universal covering, and let $f:Y\to B\pi$ be an arbitrary map. Let $w\in H^2(B\pi;\mathbb{R})$ be a cohomology class. If $c^*w\in V^2_{hyp}(X)$, then $f^*w\in V^2_{hyp}(Y)$.
\end{thm}
\begin{proof}
Without loss of generality we may assume that $c:X\hookrightarrow B\pi$ is the inclusion of a subcomplex such that the $2$-skeleton of $B\pi$ is contained in $X$. Since $X$ and $Y$ are finite there exists a finite subcomplex $X'\subset B\pi$ that contains both $X$ and $f(Y)$. We will show that $w\vert_{X'}$ is hyperbolic. Then $f^*w$ is hyperbolic by naturality and we are done.

Note that $X'$ is obtained from $X$ by attachment of finitely many cells of dimension at least $3$. By induction it suffices to consider the case where only one such cell is attached. Let $h\colon S^{k-1}\to X$ be the attaching map (with $k\geq 3$). Then $X' = X\cup_h D^k$ and $\tilde X' = \tilde X \cup_{(h\times\pi)} (D^k\times \pi)$, i.e.~the universal covering of $X'$ is obtained by attaching a $k$-cell to $\tilde X$ along each lift of $h$. 

Choose a representative $\omega\in\Omega^2(B\pi)$ of $w$. Then there is a bounded $1$-form $\alpha$ on $\tilde X$ such that $d\alpha = p^*(\omega|_X)$. Consider $\omega|_{D^k}\in\Omega^2(D^k)$. Since $H^2_{dR}(D^k)=0$ there is a $1$-form $\alpha'\in\Omega^1(D^k)$ such that $d\alpha'=\omega\vert_{D^k}$.

Now we focus on one lift of $h$ and the cell which is attached along this lift. For simplicity we will call them $h$ and $D^k$. We have $h^*\alpha\in \Omega^1(S^{k-1})$ with $d(h^*\alpha) = h^*(p^*\omega)$. Therefore, 
\[ h^*\alpha - p^*\alpha' \in \ker ( d\colon\Omega^1(S^{k-1})\to \Omega^2(S^{k-1}) ). \]
Since $k\geq 3$ the cohomology group $H^1_{dR}(S^{k-1})=0$ and there exists a function $f\in\Omega^0(S^{k-1})$ such that $df = h^*\alpha - p^*\alpha'$. Thus, $df$ is bounded by the sum of the bounds on $\alpha$ and $\alpha'$. Choose an extension of $f$ over $D^k$ that satisfies the same bound. 

We now extend $\alpha$ over $D^k$ as $p^*\alpha' + df$. Then $\alpha\in\Omega^1(\tilde X')$ is bounded and $d\alpha = p^*\omega$.
\end{proof}

The proof shows that the following extension to higher degrees is valid:
\begin{thm}\label{t:higherhominv}
Let $X$ be a finite simplicial complex such that $\pi_i(X)=0$ for $2\leq i\leq k-1$. Denote by $c\colon X\to B\pi_1(X)$ the classifying map of the universal covering. Let $f\colon Y\to B\pi$ be an arbitrary map, and let $w\in H^k(B\pi;\mathbb{R})$ be a cohomology class. If $c^*w\in V^k_{hyp}(X)$, then $f^*w\in V^k_{hyp}(Y)$.
\end{thm}

In this case we may assume that $c\colon X\to B\pi_1(X)$ is an inclusion such that the $k$-skeleton of $B\pi_1(X)$ is contained in $X$. Thus, we do not have to attach cells of dimension less than $k+1$ and the same proof as above goes through.
Note that this does not allow one to prove a formula like~\eqref{eq:hom_invar_hyp}, since it may happen that the higher skeletons of $B\pi$ can not be chosen to be finite, 
see~\cite{Stallings(1963)}.

\begin{rem}
The above Theorem~\ref{t:hominv} is very similar to Theorem~5.1 of~\cite{Kedra(2007)}. The discussion in~\cite{Kedra(2007)} is entirely in the context of manifolds, and is difficult technically. 
%for example because framings may be needed for embedded spheres whose normal bundles are not (stably) trivial. 
Our approach, extending from manifolds to simplicial complexes and formulating homological invariance in this context, is more in line with the work of the first author~\cite{Brunnbauer(2007),BrunnbauerThesis}. The generalization to Theorem~\ref{t:higherhominv} is in the same spirit as Theorem~1.9 of~\cite{BrunnbauerThesis}.
\end{rem}

%%%%%%%%%%%%%%%%%%%%%%%%%%%%%%%%%%%%%%%%%%%%%%%%%%%%%%%%%%%%%%%%%%%%%%%%%%%%%%%%%%%%%
%%%%%%%%%%%%%%%%%%%%%%%%%%%%%%%%%%%%%%%%%%%%%%%%%%%%%%%%%%%%%%%%%%%%%%%%%%%%%%%%%%%%%

\section{Amenable groups and hyperbolic classes}\label{amenable}

Consider a complete Riemannian manifold $(M,g)$. We denote by $\lambda_0(M,g)$ the largest lower bound for the spectrum of the Laplacian extended to $L^2(M)$. If $M$ is closed, then $\lambda_0(M,g)=0$ because the constant functions are in $L^2(M)$. Recall the following characterization of amenable coverings due to Brooks:
\begin{thm}[\cite{Brooks(1985)}]\label{t:amenable}
Let $(M,g)$ be a closed Riemannian manifold, and let $\bar M\to M$ be a Galois covering with Galois group $\Gamma$. Then $\Gamma$ is amenable if and only if $\lambda_0(\bar M,\bar g)=0$, where $\bar g$ denotes the lifted metric.
\end{thm}

We now use this result to prove the following:
\begin{thm}\label{amenable->hyp=0}
Let $X$ be a finite simplicial complex with amenable fundamental group. For all $k$ we have 
\[ V^k_{hyp}(X)=0 \ . \]
\end{thm}
\begin{proof}
Assume there is a nontrivial $w\in V^k_{hyp}(X)$. By the well known result of Thom, there is a map $f\colon N\to X$ from a connected closed orientable $k$-dimensional manifold $N$ such that $f^*w\neq 0\in H^k(N;\mathbb{R})$. We may assume without loss of generality that $f_*\colon\pi_1(N)\to\pi_1(X)$ is surjective. 

Consider the Galois covering $\bar N= f^*\tilde X$ of $N$. Its Galois group is $\pi_1(X)$, which by assumption is amenable. Thus, by Theorem~\ref{t:amenable}, we have $\lambda_0(\bar N, \bar g)=0$ for any metric $g$ on $N$.

Recall that the isoperimetric constant of a complete $k$-dimensional manifold $(\bar N, \bar g)$ is defined as
\[ i(\bar N, \bar g)= \inf_\Omega \frac{\Vol_{k-1}(\partial\Omega)}{\Vol_k(\Omega)}, \]
where the infimum is taken over all relatively compact sets $\Omega\subset\bar N$ with sufficiently regular boundary. The Cheeger inequality~\cite{Cheeger(1970)} tells us that 
\[
{\textstyle\frac{1}{4}}i(\bar N, \bar g)^2\leq \lambda_0(\bar N, \bar g) \ ,
\]
so that we conclude $i(\bar N, \bar g)=0$. This will lead to a contradiction.

Suppose that $f^*w$ is represented by the volume form $\omega$ of a Riemannian metric $g$ -- this is possible because $f^*w\neq 0$ in the top-degree cohomology of $N$. Then, on the one hand, $p^*\omega$ is the Riemannian volume form of $\bar g$. On the other hand, $p^*\omega = d\alpha$ for some $\alpha$ which is bounded, as we see from the commutativity of the following diagram:
\[
\xymatrix{
\bar N \ar_-p[d] \ar[r] & \tilde X \ar^-p[d] \\
N \ar^-f[r] & X 
}
\]
Now for any $\Omega\subset\bar N$ with sufficiently smooth boundary we can apply Stokes's theorem to obtain
\[
\Vol_k(\Omega) = \int_\Omega p^*\omega = \int_{\partial\Omega}\alpha\leq c \Vol_{k-1}(\partial\Omega) \ ,
\]
where $c$ is any $C^0$-bound for $\alpha$. It follows that $i(\bar N, \bar g)\geq 1/c$, contradicting the vanishing of the isoperimetric constant. 
This contradiction completes the proof.
\end{proof}

\begin{rem}
Proposition~1.9 and Theorem~6.7 of~\cite{Kedra(2007)} are special cases of the above Theorem~\ref{amenable->hyp=0}.
\end{rem}

The converse of Theorem~\ref{amenable->hyp=0} is not true for finite simplicial complexes. For example, if $M$ is a closed hyperbolic three-manifold which is a real homology sphere, and $X$ is obtained from $M$ by removing the interior of a top-degree simplex, then $V^k_{hyp}(X)=0$ for all $k$, but the fundamental group of $X$ is non-amenable. However, if we restrict to closed oriented manifolds, then there is the following strong converse of Theorem~\ref{amenable->hyp=0} due to Gromov, Brooks, Sikorav and others:
\begin{thm}\label{t:conv}
Let $M$ be a closed oriented smooth $n$-manifold. If $V^n_{hyp}(M)=0$, then $\pi_1(M)$ is amenable.
\end{thm}
\begin{proof}
One has to prove that if $\pi_1(M)$ is not amenable, then the lift of the volume form of any Riemannian metric $g$ on $M$ to $\tilde M$ admits a bounded primitive. By Theorem~\ref{t:amenable}, the non-amenability of $\pi_1(M)$ is equivalent to the non-vanishing of $\lambda_0(\tilde M,\tilde g)$. Using the converse of the Cheeger inequality due to Buser~\cite{Buser(1982)}, this implies the non-vanishing of the isoperimetric constant $i(\tilde M,\tilde g)$. As explained in~\cite{Brooks(1983)}, \cite[Chapter 6]{Gromov(1999)} and~\cite{Sik01}, the non-vanishing of $i(\tilde M,\tilde g)$ leads to the existence of a bounded primitive for the volume form. See the paper of Sikorav~\cite{Sik01} for a detailed proof not passing through Theorem~\ref{t:amenable}, and  Block and Weinberger~\cite{BW92} for a different approach.
\end{proof}

\section{Examples and applications}\label{s:examples}

Gromov~\cite[Section 6C]{Gromov(1993)} showed that the inclusion $V^k_b(X) \subset V^k_{hyp}(X)$ is usually strict in degrees $k\geq 3$. His examples are of the following form: choose a closed orientable manifold $M$ of dimension $k-1$ with non-amenable fundamental group (this is where $k\geq 3$ is used), and take $X$ to be the product $M\times S^1$. The fundamental group of $X$ is non-amenable and therefore the top degree cohomology is equal to its hyperbolic subspace by Theorem~\ref{t:conv}. But since the simplicial volume $\|X\|$ is zero due to the  presence of a free circle action, it follows that $V^k_b(X)=0$ (see \cite{Gromov(1983)}, page 17). This phenomenon shows in particular that Theorem~\ref{amenable->hyp=0} is not a consequence of the vanishing theorem for the bounded cohomology of amenable groups.

It is an open question whether the inclusion $V^2_b(X) \subset V^2_{hyp}(X)$ may also be strict. In fact, Gromov~\cite[Section 6C]{Gromov(1993)} conjectured that this inclusion is never strict, which, together with the above examples for higher degrees, would completely resolve Question~1.15 of~\cite{Kedra(2007)}.  

K\c{e}dra~\cite[Question 1.10]{Kedra(2007)} also asked whether every atoroidal cohomology class of degree two is hyperbolic. We now give a negative answer to this question using the following result of Barge and Ghys:
\begin{thm}[\cite{BG(1988)}]
For every positive integer $k$ there exists a finitely presentable nilpotent group $\Gamma$ such that $H_2(B\Gamma;\mathbb{Z})$ contains a non-torsion element $a$ that is not contained in the subgroup generated by all elements which are representable by surfaces of genus at most $k$.
\end{thm}

Note that nilpotent groups are amenable. Therefore, for these groups $V^2_{hyp}(B\Gamma) = 0$ by Theorem~\ref{amenable->hyp=0}. Consider a class $w\in H^2(B\Gamma;\mathbb{R}) = \Hom (H_2(B\Gamma;\mathbb{Z}),\mathbb{R})$ that sends the subgroup generated by all elements which are representable by surfaces of genus at most $k$ to zero but that fulfills $w(a)\neq 0$. If $k\geq 1$, then $w$ is atoroidal. Thus, we have:
\begin{cor}
There exist finitely presentable groups $\Gamma$ such that $V^2_{hyp}(B\Gamma)=0$ and $V^2_{ator}(B\Gamma)\neq 0$.
\end{cor}

These examples can be realized by symplectic forms on closed manifolds:

\begin{cor}\label{c:sympl}
There exist closed symplectic four-manifolds $(M,\omega)$ for which (the cohomology class of) $\omega$ is atoroidal but not hyperbolic.
\end{cor}
\begin{proof}
We use the same construction as in Section~3.3 of~\cite{Kedra(2007)}, compare also~\cite{ABKP00,KRT(2007)}.

Let $\Gamma$ be one of the groups from the construction of Barge and Ghys~\cite{BG(1988)}, and $w\in H^2(B\Gamma;\mathbb{R})$ a non-torsion atoroidal class. By the construction of Amor\'os et.~al.~\cite{ABKP00} there exists a Lefschetz fibration $N$ over $S^2$ with a section, with $\pi_1(N)=\Gamma$, and such that $\omega_F:=c^*w$ evaluates non-trivially on the fiber $F$ of $N\to S^2$. Let $B$ be a surface of genus at least $2$, and $M$ the fiber sum of $N$ with $B\times F$ along $F$. Then $\pi\colon M\to B$ is a Lefschetz fibration with $\pi_1(M)=\Gamma\times\pi_1(B)$. If $\omega_B$ denotes a generator for $H^2(B;\mathbb{R})$, then for all $k\in\mathbb{R}$ which are large enough, the cohomology class $\omega:=\omega_F+k\pi^*\omega_B$ is represented by a symplectic form constructed by the generalization of the Thurston construction from bundles to Lefschetz fibrations. On the one hand, $\omega_F$ and $\omega_B$ are both atoroidal, and therefore so is their linear combination $\omega$. On the other hand, if $\omega$ were hyperbolic, then, because $\omega_B$ is hyperbolic, it would follow that $\omega_F$ would be hyperbolic, which would be a contradiction.
\end{proof}

\begin{rem}
In the examples constructed in the proof, more is true than was claimed in the statement of Corollary~\ref{c:sympl}. Namely, not only is the class $[\omega]$ not hyperbolic, but the cohomology class of any symplectic form on a manifold $M$ with the given fundamental group $\Gamma\times\pi_1(B)$ fails to be hyperbolic. This is because $V^2_{hyp}(M)=\pi^*V^2_{hyp}(B)$ is an isotropic subspace for the cup product on $H^2(M;\mathbb{R})$, and so can not contain the class of any non-degenerate two-form. (Here $\pi$ is the composition of the classifying map of $M$ with the projection $B\Gamma\times B\to B$.) 
\end{rem}

\begin{rem}
The examples in Corollary~\ref{c:sympl} exhibit the same behavior as Gromov's examples mentioned above: in top degree the hyperbolic subspace $V^4_{hyp}(M)$ is non-trivial, but the bounded
subspace $V^4_{b}(M)$ vanishes. The former is due to the hyperbolicity of the volume form $\omega_F\wedge\pi^*\omega_B$, which follows from the hyperbolicity of $\omega_B$, without even
appealing to Theorem~\ref{t:conv}.
The latter is due to the amenability of $\Gamma$, which implies the vanishing of its bounded cohomology. Note that the sum of all the hyperbolic subspaces is always an ideal in the real cohomology ring,
but this is not true for the bounded subspaces.
\end{rem}

\bigskip

%%%%%%%%%%%%%%%%%%%%%%%%%%%%%%%%%%%%%%%%%%%%%%%%%%%%%%%%%%%%%%%%%%%%%%%%%%%%%%%%%%%%%
%%%%%%%%%%%%%%%%%%%%%%%%%%%%%%%%%%%%%%%%%%%%%%%%%%%%%%%%%%%%%%%%%%%%%%%%%%%%%%%%%%%%%

%\bibliographystyle{amsalpha}
%\bibliography{asymp_invar}  

%\providecommand{\bysame}{\leavevmode\hbox to3em{\hrulefill}\thinspace}
\providecommand{\MR}{\relax\ifhmode\unskip\space\fi MR }
\providecommand{\MRhref}[2]{%
  \href{http://www.ams.org/mathscinet-getitem?mr=#1}{#2}
}
\providecommand{\href}[2]{#2}

\end{document}